\documentclass{amsart}
\usepackage{amscd,amsthm}
\usepackage{latexsym}
\usepackage{capt-of}
\usepackage{graphicx}
\usepackage{algorithm}
\usepackage[noend]{algpseudocode}

\makeatletter
\def\BState{\State\hskip-\ALG@thistlm}
\makeatother

\DeclareFontFamily{U}{wncy}{}
\DeclareFontShape{U}{wncy}{m}{n}{<->wncyr10}{}
\DeclareSymbolFont{mcy}{U}{wncy}{m}{n}
\DeclareMathSymbol{\Sha}{\mathord}{mcy}{"58}

\usepackage{amssymb,amsmath}
\usepackage{amsfonts}
\usepackage{hyperref}
\newcounter{ctfig}

\newcommand{\C}{\mathcal{C}}

\newcommand{\F}{\mathcal{F}}

\newcommand{\JacC}{{\hbox{Jac}_{\lower.5pt\hbox{$_\C$}}}}
\newcommand{\JacF}{{\hbox{Jac}_{\lower.5pt\hbox{$_\F$}}}}

\theoremstyle{plain}

\theoremstyle{definition}

\def\F{{\mathbb F}}

\def\C{{\mathbb C}}

\def\e32{{{}_3E_2}}
\def\f32{{{}_3F_2}}
\def\a32{{{}_3A_2}}

\begin{document}
%% \foreach \x in{graphics,floats}{%
%%     \immediate\write18{pdflatex -jobname=template-\x\space "\def\noexpand\placeholder{\x} \noexpand\input{template}"}%
%%     \includepdf[pages=-]{template-\x}%
%% }
\bibliographystyle{plain}
\bibstyle{plain}

\title[Sums of Three Cubes]{Representing integers as a sum of three cubes}

\author{Jon Grantham}
\email{grantham@super.org}
\address{IDA-CCS 17100 Science Dr. Bowie, MD 20715 USA}

\author{P.G. Walsh}
\email{gwalsh@uottawa.ca}
\address{Department of Mathematics and Statistics, University of Ottawa}
\subjclass[2020]{11D25}
\keywords{diophantine equations}
\date{\today}

\begin{abstract}
In this article we further develop methods for representing integers as a sum of three cubes. In particular, a barrier to solving the case $k=3$, which was outlined in a previous paper of the second author, is overcome. A very recent computation indicates that the method is quite favourable to other methods in terms of time estimates. A hybrid of the method presented here and those in a previous paper is currently underway for unsolved cases.
\end{abstract}

\maketitle

\section{Introduction}
The problem of representing an integer as a sum of three cubes has drawn much attention in recent years with the recent discoveries of Booker and Sutherland by way of solving certain notoriously famous cases going back to a problem posed by Mordell in the early 1950's. For more on the history of the problem, the reader is referred to \cite{W} and \cite{W1}.\\

The approach taken in \cite{W} to solve the equation
$$X^3+Y^3+Z^3=k \leqno (1.1)$$
made use of the arithmetic in the field $\mathbb{Q}(k^{1/3})$ in order to search for the quantity $X+Y$. If the correct value for $X+Y$ had been found, it was shown that a relatively small computation would then find $Z$ in a certain arithmetic progression, from which $X$ and $Y$ could then be computed very easily. As $X+Y$ is a divisor of $(-Z)^3+k$, the search for the correct value of $X+Y$ was taken over integers which are norms of elements 
$$\gamma = a+b\cdot k^{1/3}+c\cdot k^{2/3} \in \mathbb{Z}[k^{1/3}],$$
or integers of that form divided by very small factors. In particular, it was shown in \cite{W} how a solution $X,Y,Z$ to the original problem arose from such algebraic numbers $\gamma$, together with another integer $s$, which located $Z$ inside of an arithmetic progression, and that for many of the published very large solutions of (1.1), the height of the associated $\gamma$, and value of $s$, are strikingly small.\\

Two strategic ideas were presented in \cite{W} to assist in searching for the correct value of $X+Y$. The first, alluded to above, is to deal with the case that the class number of the cubic field generated by $k^{1/3}$ is larger than $1$. The second strategic idea served the purpose of reducing the set of possible values $a,b,c$ by exploiting a dependency of one of the three of these values on the other two. This idea was shown to be very useful for certain known large solutions to cases such as $k=33$. However, it was shown in \cite{W} that in the case $k=3$, not only was this idea shown to be of no help, the further problem that the value $s$ described above is quite large, adding a considerable component to the overall running time.\\

The purpose of the paper is to develop a new strategy. In so doing, we will see that we can circumvent the problems discussed above, and in fact arrive at a run time which is very competitive with that in the work of Booker and Sutherland. In their paper, it is estimated that to solve the case $k=3$, several hundred core years are required. We estimate that the algorithm we present in this paper should take about $137$ core years.

\section{Algorithm Overview and an Example}
In what follows, $k$ will represent a positive non-cube integer, and let $\alpha = k^{1/3}$. We first describe the approach in the case that the field $K=\mathbb{Q}(\alpha )$ has class number one. For certain reasons that are beneficial to us, we rewrite $Y$ and $Z$ in (1.1) as $-Y$ and $-Z$ so that (1.1) can conveniently be rewritten as
$$X^3-Y^3=Z^3+k. \leqno (2.1)$$
We define the integer $t$ to represent the difference $X-Y$. With the hypothesis made on the class number, it follows that there are integers $a,b,c$ for which
$$t=N_{K/\mathbb{Q}}(a+b\alpha +c\alpha^2).$$
Since $X=Y+t$, it also follows that there are integers $d,e,f$ for which
$$X^2+XY+Y^2=3Y^2+3Yt+t^2=N_{K/\mathbb{Q}}(d + e\alpha + f \alpha^2). \leqno (2.2)$$
The key now is to write $d,e,f$ in terms of $a,b,c$ which can be done using the equation
$$(a+b\alpha +c\alpha^2)(d + e\alpha + f \alpha^2)=Z+\alpha,$$
along with a little bit of linear algebra. In fact, it is quite easy to verify that\newline
$\; $\newline
\centerline{$d=\frac{1}{t}(Z(a^2-kbc)+(kb^2-kac)),$}\newline
\centerline{$e=\frac{1}{t}(Z(kc^2-ac)+(a^2-kbc)),$}\newline
\centerline{$f=\frac{1}{t}(Z(b^2-ac)+(kc^2-ab)).$}\newline
$\; $\newline
If it's the case that $t \ge 2$, $\gcd(t,a^2-kbc) =1$, and $Z=z_1$ is defined to be
$$ z_1 = -(kb^2-kac)(a^2-kbc)^{-1} \; (\bmod \; t),$$
then $d$, $e$ and $f$ are integers. In particular, as $Z=z_1+st$ for some (unknown) integer $s$, $d$, $e$, and $f$ become linear polynomials in $s$ with integer coefficients, and the equality on the right in (2.2) is an integral quadratic polynomial in $Y$ equaling an integral cubic polynomial in $s$. Taking this just one small step further, we see that
$$(6Y+3t)^2=12N_{K/\mathbb{Q}}(d + e\alpha + f \alpha^2)-3t^2. \leqno (2.3)$$
Therefore, $(s,6Y+3t)$ is an integral point on an elliptic curve which is determined entirely by $a,b,c$.\\

\noindent {\bf Example 2.1.} For further clarity, let us examine $k=3$. In this case, the relevant field $\mathbb{Q}(3^{1/3})$ has class number one, and so the pertinent remarks above apply. The solution $(X,Y,Z)$ found by Booker and Sutherland is given by $(X,Y,Z)$
$$=(569936821221962380720,569936821113563493509,472715493453327032).$$
The quantity $t=X-Y$ associated to this solution is 
$t=108398887211$, which the reader will notice is relatively small, as $X$ and $Y$ are remarkably close together for their size. Our algorithm will find $t$ once a triple $a,b,c$ is found for which $t=N_{K/\mathbb{Q}}(a+b\alpha +c\alpha^2).$ It so happens that there are several such triples with all of $|a|,|b|,|c|$ relatively small. In fact, three such triples are given by
$$\{ (5603, 1612, 1156), (-2386, -143, 2432), (3824, -2135, 992)\}.$$
We will focus on the first of these. With $(a,b,c)=(5603,1612,1156)$, we first determine that $\gcd(t,a^2-3bc)=\gcd(108398887211,25803193)=1$, and so $z_1$ is computed to be $-3(b^2-ac)(a^2-3bc)^{-1} \; (\bmod \; t) = 87001523664$, from which the linear functions $d,e,f$ are determined:\newline
$\; $\newline
\centerline{$d=25803193s+20709780,$}\newline
\centerline{$e=-5023028s-4031509,$}\newline
\centerline{$f=-3878524s-3112924.$}\newline
$\; $\newline
Using (2.2) and solving the square on the left hand side gives that $(s,6Y+3t)$ is an integral point on the curve\newline
\centerline{$r^2=141003824982997192302252s^3$}\newline
\centerline{$+ 339511260630207004599744s^2$}\newline
\centerline{$+ 272493544314871871456256s$}\newline 
\centerline{$+ 37650633287640319405761.$}\\
$\; $\newline
The input value $s=4360888$ satisfies the property that the value of the polynomial is a square, with $y$-coordinate having the property that upon subtracting $3t$, it is divisible by $6$. Thus, we finally obtain $Z=z_1+st=472715493453327032$, $Y=569936821113563493509$, and $X=Y+t$.\\

With this method, the pari-gp routine {\em hyperellratpoints} appears to be an extremely promising way to approach the problem. The throughput rate on the second author's pc, which is equipped with an i5 processor and 32 GB of ram, is roughly $3$ milliseconds per curve on one core of the processor. With outer loop being $c$, the coefficient of $3^{2/3}$, the total number of curves is roughly $10^{11}$, which all translates into about an 11 core-year computation, or almost a two-year computation on this single pc. This compares favourably with the computing estimate given in \cite{BS}.\\

\section{Class Number Considerations}

As discussed in \cite{W}, it can often be the case that the field in question does not have unique factorization, and in that case, the value $t=X-Y$ we are searching for will not show up as a norm. However, under GRH, as discussed in \cite{Sar}, a small multiple of $t$ is likely to be a norm of an element $\gamma = a+b\alpha + c\alpha^2$. Thus, to circumvent the problem of non-unique factorization, one can insert a short loop over small integers $l$ so that if $g=\gcd(t,l)$ and $N_{K/\mathbb{Q}}(\gamma) = t$, then $X=Y+(t/g)$, $Y$, $Z=z_1+st$ will be a solution to (2.1) provided
$$(6Y+3(t/g))^2 = 12gN_{K/\mathbb{Q}}(d+e\alpha+f\alpha^2)-3(t/g)^2,$$
where all of the quantities above are defined in Section 2.\\

The reader is invited to access the pari code available at \cite{W2}, modify the parameters therein, and perform these computations. In our experience, all of the moderately sized solutions to (1.1) are found quite quickly. The reader is also referred to the website \cite{H} containing a catalogue of solutions to (1.1) with $\max (|X|,|Y|,|Z|) < 10^{15}$, which is helpful when comparing results of computation using the method presented here. We issue a warning however that the monumental task of producing the enormous table in \cite{H} has left open only very large problems, which even our method will require a significant amount of running time to expand upon. For example, our significant effort to solve the case $k=114$ remains unrewarded.\\

\bibliography{art}

\end{document}